%% file: bettersieve.tex
\documentstyle[12pt, fullpage]{amsart}
\begin{document}
\input{universal}
\createtitle{The Least Prime Primitive Root and the Shifted
Sieve}{11N69, 11N36}

\def\w{{\omega(\phi(q))}}
\def\f{{\phi(q)}}
\def\ind#1{{\mathop{\rm ind}#1}}
\def\nup{{{\nu\!}_p}} \def\nud{{{\nu\!}_d}}
\def\Z#1{{\bf Z}_{#1}}
\def\Zx#1{{\bf Z}^\times_{#1}}
\def\lr{$\lambda$\discretionary{-}{}{-}root}
\def\lamm#1{{\lambda_{#1}^-}}

\def\Q#1#2{{\cal Q}(#1,#2)}
\def\Qs#1{\Q\sigma{#1}}
\def\Qp#1#2{{\cal Q}'(#1,#2)}
\def\Qps#1{\Qp\sigma{#1}}
\def\R#1#2{{\cal R}(#1,#2)}
\def\Bs{{4/(1-\sigma)}}
\def\A{{\cal A}}

\def\psip#1#2{{\psi_1(#1,#2)}}
\def\starsum{\mathop{{\sum}^*}}

\section{Introduction}\noindent
\input{intro}
\section{Preliminaries}\label{charsec}\noindent
\input{finite}

\input{charsum}
\section{The shifted sieve: Proof of Proposition
 \ref{gqprop}}\label{Rossec}\noindent
\input{shifted}
\input{applyshift}
\section{Proof of Lemmas \ref{aaplem} and \ref{aalem}}\label{aasec}\noindent
\input{Zhang}

\input{almostall}

\bibliography{bettersieve}
\bibliographystyle{amsplain}
\end{document}

%% file: universal.tex
\newcommand{\createtitle}[2]{\title{#1}\author{Greg Martin}\address{Department of Mathematics\\University of Toronto\\Canada M5S 3G3}\email{gerg@@math.toronto.edu}\subjclass{#2}\maketitle}
\newcommand{\make}[1]{\label{#1sec}\noindent\input{#1}}
\newcommand{\Mmake}[1]{\label{#1sec}\noindent}

%%% use for AMSTeX
% \input{/opt/texmf/tex/ams/amstex}
%%%

\newtheorem{theorem}{Theorem}
\newtheorem{lemma}[theorem]{Lemma}
\newtheorem{corollary}{Corollary}[theorem]
\newtheorem{proposition}[theorem]{Proposition}

\newenvironment{pflike}[1]{\noindent{\bf #1}}{\vskip10pt} % NO \QED
\newenvironment{proof}{\begin{pflike}{Proof:}}{\qed\end{pflike}}

\newcommand{\2}[1]{\ifmmode{\cal#1}\else$\cal#1$\fi}
\newcommand{\3}[1]{\#\{#1\}}
\newcommand{\abs}[1]{\left|#1\right|} % take out
\newcommand{\floor}[1]{\lfloor#1\rfloor}
\newcommand{\bfloor}[1]{\big\lfloor#1\big\rfloor}
\newcommand{\bbfloor}[1]{\bigg\lfloor#1\bigg\rfloor}
\newcommand{\ceil}[1]{\lceil#1\rceil}
\newcommand{\bceil}[1]{\big\lceil#1\big\rceil}
\newcommand{\bbceil}[1]{\bigg\lceil#1\bigg\rceil}

\renewcommand{\mod}[1]{{\ifmmode\text{\rm\ (mod~$#1$)}
 \else\discretionary{}{}{\hbox{ }}\rm(mod~$#1$)\fi}}
\newcommand{\ep}{\varepsilon}
\renewcommand{\implies}{\Rightarrow}
\newcommand{\rmif}{{\rm if\ }}

\newcommand{\half}{{\mathchoice{\textstyle\frac12}{1/2}{1/2}{1/2}}}
\newsymbol\dnd 232D % see The Joy of TeX (2nd ed.), Appendix G
\newcommand{\exdiv}{\mathrel{\mid\mid}}

\renewcommand{\lg}[1]{\mathop{\log_{#1}}}
\def\lgs#1^#2{\mathop{\log_{#1}^{#2}}}
\newcommand{\li}{\mathop{\rm li}}

\newcommand{\doublespace}{%\parskip 0.2cm
  \baselineskip=24pt}
\newcommand{\spaceandahalf}{\parskip6pt\baselineskip=18pt}

\newcommand{\scroll}[1]{\scrollmode#1\errorstopmode}
\newcommand{\comment}[1]{}

\vfuzz=2pt % prohibits overfull \vbox messages due to page headers
	   % and footers

%% file: intro.tex
If $p$ is a prime, we define $g^*(p)$ to be the least prime that is a
primitive root\mod{p}, and similarly for prime powers $p^r$. The
problem of establishing a bound for $g^*(p)$ uniformly in $p$ is quite
difficult, comparable with establishing a uniform upper bound for the
least prime in an arithmetic progression. Indeed, there do not exist
any uniform upper bounds for $g^*(p)$ that improve upon the current
bounds for the least prime in an arithmetic progression. However, much
more can be said if we exclude a very small set of primes. The purpose
of this paper is to improve existing bounds for $g^*(p)$ which hold
for almost all primes $p$, and to establish analogous results for all
composite moduli.

Elliott \cite{Ell:TDoPR} had first given a bound for $g^*(p)$ for all
but $O(Y^\ep)$ primes $p$ up to $Y$, of the form $g^*(p)\le(\log
p)^{O_\ep(\lg3p)}$. (Here we have defined $\lg1x=\max\{\log x,1\}$ and
$\lg nx=\max\{\log(\lg{n-1}x),1\}$ for any integer $n\ge2$.) This was
subsequently improved by Nongkynrih \cite{Non:OPPR} to $g^*(p)\le(\log
p)^{O_\ep(\lg3p/\lg4p)}$. We are able to establish the following
bound. Write $\omega(n)$ for the number of distinct prime factors of
$n$.
\begin{theorem}
Let $Y$, $\ep$, and $\eta$ be positive real numbers with
$\ep\le20/21$, and define $B=B(\ep,\eta)=\frac3\ep+\frac54+\eta$. The
number of odd prime powers $p^r$ not exceeding $Y$ for which the
estimate
\begin{equation*}
g^*(p^r) \ll_{\ep,\eta} \big( \omega(p-1)^2\log p \big)^B
\end{equation*}
fails is $O_{\ep,\eta}\big(Y^\ep\big)$.
\label{aapthm}
\end{theorem}
Since $\omega(n)\ll\log n$ for all integers $n$, it is apparent that
the bound for $g^*(p^r)$ given in Theorem \ref{aapthm} is no larger
than a fixed (depending on $\ep$ and $\eta$) power of $\log p$. We see
that this is an improvement over the existing bounds, where the
exponent of $\log p$ tends to infinity with $p$. We remark that
Theorem \ref{aapthm} may easily be extended to include all moduli
which admit primitive roots, i.e., to include moduli of the form
$2p^r$.

To extend this type of result to composite moduli, we use the
following definition. Given an integer $q\ge2$, we say that a
{\it\lr\mod q} is an integer, coprime to $q$, whose multiplicative
order is maximal among all integers coprime to $q$. We see that the
\lr{} is an extension of the primitive root to all moduli, and we
extend the notation $g^*(q)$ to mean the least prime \lr\mod q.
\begin{theorem}
Let $\ep$ be a positive real number. For almost all integers $q\ge2$,
we have
\begin{equation*}
g^*(q) \ll_\ep \omega(\phi(q))^{44/5+\ep}(\log q)^{22/5}.
\end{equation*}
\label{aathm}
\end{theorem}
The approach to establishing these theorems is through Proposition
\ref{gqprop} below, which gives a bound for $g^*(q)$ based on the
assumption of a zero-free rectangle for characters\mod q. This is the
same approach taken in earlier work on this subject; the improvement
lies in the use of the ``shifted sieve'', a version of the linear
sieve with very good error terms, rather than Brun's sieve.

For any integer $n$, let $s(n)$ denote the largest squarefree divisor
of $n$. For any integer $q\ge2$, let $E(q)$ denote the exponent of the
group $\Zx q$ of reduced residue classes\mod q, let $\Phi(q)$ be the
group of Dirichlet characters\mod q, and define
\begin{equation*}
\Phi_*(q) = \{\chi^{E(q)/s(\f)}: \chi\in\Phi(q)\}.
\end{equation*}
Only the characters in $\Phi_*(q)$ are relevant to detecting
\lr{s}, as we show in Section \ref{charsec}. Let $c_0$ be the
probability that a randomly chosen element of $\Zx q$ is a \lr. Also,
given real numbers $\sigma$ and $T$ with $\half\le\sigma<1$ and $T>0$,
define $\Qs T$ to be the set of integers $q\ge2$ such that, for some
nonprincipal $\chi\in\Phi_*(q)$, the corresponding $L$-function
$L(s,\chi)$ has a zero $\beta+i\gamma$ with $\beta>\sigma$ and
$\abs\gamma<T$.
\begin{proposition}
Let $q\ge2$ be an integer and $\sigma$ a real number satisfying
$\half\le\sigma<1$, and set
\begin{equation*}
f(q,\sigma)= \big( \w^2\lg1\w \cdot c_0^{-1}\log q \big)^{1/(1-\sigma)}.
\end{equation*}
If $q\notin\Qs{f(q,\sigma)}$, then $g^*(q)\ll_\sigma f(q,\sigma)$.
\label{gqprop}
\end{proposition}
We remark that $f(q,\sigma)\ll_{\sigma,\theta} q^\theta$ for every
$\theta>0$. We also remark that $c_0^{-1}\ll\lg1\w$ (see Section
\ref{charsec}) and that the generalized Riemann hypothesis implies
that $\Q\half T$ is empty for every $T>0$. Thus the following
corollary of Proposition \ref{gqprop} is immediate.
\begin{corollary}
If the generalized Riemann hypothesis holds for (certain)
characters\mod q, then
\begin{equation*}
g^*(q)\ll \big(\w\lg1\w\big)^4(\log q)^2.
\end{equation*}
\end{corollary}
In the case where $q$ is a prime, this has already been shown by Shoup
\cite{Sho:SfPRiFF}, improving an earlier result of Wang
\cite{Wan:OtLPRoaP} in which $(\w\lg1\w)^4$ is replaced by
$\w^6$. Although both authors state their bounds only for primitive
roots, the bounds actually hold for prime primitive roots as well.

To deduce Theorems \ref{aapthm} and \ref{aathm} from Proposition
\ref{gqprop}, we need bounds on the size of $\Qs T$. To this end, we
define $Q(Y;\sigma,T)$ to be the number of elements of $\Qs T$ not
exceeding $Y$, and $Q'(Y;\sigma,T)$ to be the number of elements of
$\Qs T$ which are odd prime powers not exceeding $Y$. The following
lemmas, when combined with Proposition \ref{gqprop}, imply Theorems
\ref{aapthm} and \ref{aathm}.
\begin{lemma}
Let $Y$, $\ep$, $\eta$, and $B$ be as in Theorem \ref{aapthm}. There
exists $\theta=\theta(\ep,\eta)>0$ such that
\begin{equation*}
Q'(Y; 1-B^{-1}, Y^\theta) \ll_{\ep,\eta} Y^\ep.
\end{equation*}
\label{aaplem}
\end{lemma}
\vskip-12pt
\begin{lemma}
We have $Q(Y;{17\over22},Y^{1/20})=o(Y)$.
\label{aalem}
\end{lemma}
Lemma \ref{aaplem} follows directly from existing zero-density
estimates for Dirichlet $L$-functions, but Lemma \ref{aalem} is
somewhat more complicated due to the prevalence of imprimitive
characters in $\Phi_*(q)$ for composite moduli $q$ (see Section
\ref{aasec}).

The author would like to express his gratitude to Hugh Montgomery for
suggesting this problem and to thank him and Trevor Wooley for their
guidance and support. The author would also like to thank Andrew
Granville and Andrew Odlyzko for their comments regarding existing
results related to this work. This material is based upon work
supported under a National Science Foundation Graduate Research
Fellowship.

%% file: finite.tex
We begin by developing some notation and simple facts relating to the
characters\mod q which are relevant to detecting \lr{s}. Let $G$ be a
finite abelian group with exponent $E$. For every prime $\ell$ that
divides $E$, let $\alpha(\ell)$ be the largest integer such that
$\ell^{\alpha(\ell)}$ divides $E$. There exist integers $m(\ell)$ for
which we can write
\begin{equation*}
G \cong \bigg( \bigoplus_{\ell\mid E} \big( \Z{\ell^{\alpha(\ell)}}
\big)^{m(\ell)} \bigg) \bigoplus H
\end{equation*}
for some subgroup $H$ whose exponent divides $E/s(E)$. For each prime
$p$ dividing $E$, we define subgroups $G_p$ of $G$ by
\begin{equation}
G_p = \big( p\Z{p^{\alpha(p)}} \big)^{m(p)} \bigoplus \bigg( \bigoplus
\begin{Sb}\ell\mid E \\ \ell\ne p\end{Sb} \big( \Z{\ell^{\alpha(\ell)}}
\big)^{m(\ell)} \bigg) \bigoplus H,
\label{Gpdef}
\end{equation}
the set of all elements of $G$ whose order divides $E/p$. We see that
the index of $G_p$ in $G$ is $p^{m(p)}$. We extend this notation to
all squarefree divisors $d$ of $E$ by defining subgroups $G_d$ by
\begin{equation*}
G_d = \bigcap_{p\mid d} G_p,
\end{equation*}
and (abusing notation somewhat) we define $m(d)$ to be the real number
which satisfies
\begin{equation*}
d^{m(d)} = \prod_{p\mid d} p^{m(p)},
\end{equation*}
so that $d^{m(d)}$ is a multiplicative function of $d$. By convention,
we let $G_1=G$ and $m(1)=1$. We note that $m(d)\ge1$ for all
squarefree divisors $d$ of $E$, and that the index of $G_d$ in $G$ is
$d^{m(d)}$.

Let $\gamma(g)$ be the characteristic function of elements of maximal
order in $G$. Then, by definition (\ref{Gpdef}) of the $G_p$, we have
\begin{equation}
\{ g\in G: \gamma(g)=1 \} = G \setminus \bigcup_{p\mid E}
G_p.  \label{Gexclude}
\end{equation}
If we define $\nu(g)$ to be the product of all primes $p$ dividing $E$
such that $g\in G_p$ (or equivalently, the largest squarefree divisor
$d$ of $E$ such that $g\in G_d$), then we see from equation
(\ref{Gexclude}) that for any $g\in G$, we have
\begin{equation}
\gamma(g) = \begin{cases} 1&\rmif \nu(g)=1, \\ 0&\rmif \nu(g)>1.
\end{cases}  \label{manyBiggs2}
\end{equation}

We may also detect these elements of maximal order using group
characters. Let $\Phi$ be the group of homomorphisms from $G$ into
$\bf C$. For each squarefree $d$ dividing $E$, define subgroups
$\Phi_d$ of the character group $\Phi$ by
\begin{equation*}
\Phi_d = \big\{ \chi^{E/d}: \chi\in\Phi \big\}. \label{Gddef}
\end{equation*}
For convenience we write $\Phi_*$ for $\Phi_{s(E)}$. Let $h_d$ be the
characteristic function of $G_d$. By the standard properties of group
characters, for any $g\in G$ we have
\begin{equation}
h_d(g) = \frac1{\abs{\Phi_d}} \sum_{\chi\in\Phi_d} \chi(g).
\label{charfnGd}
\end{equation}
By summing this over all $g\in G$ we see that $|\Phi_d|=\abs
G/\abs{G_d}=d^{m(d)}$, and in fact we can treat this as the definition
of the real numbers $m(d)$. Finally, we define $c_0$ to be the
probability that a randomly chosen element of $\Zx q$ is a \lr. From
equation (\ref{Gexclude}) and the definition (\ref{Gpdef}) of the
$G_p$, we can easily calculate that
\begin{equation*}
c_0 = \prod_{p\mid\phi(q)} \big( 1-\frac1{p^{m(p)}} \big).
\end{equation*}
We note in particular that $c_0^{-1}\le q/\phi(q)\ll\lg1\w$.

%% file: charsum.tex
In the course of applying the sieve, it will be important to
understand the behavior of the sum $\psi_1(x,\chi)$ defined by
\begin{equation*}
\psi_1(x,\chi) = \sum_{n<x} \chi(n)\Lambda(n)(x-n).
\end{equation*}
The following lemma provides the necessary bound, for the moduli $q$
for which Proposition \ref{gqprop} will be established.
%Let $\psip x{\chi_0}=\psi_1(x,\chi_0)-\half x^2$, while for
%$\chi\ne\chi_0$, let $\psip x\chi=\psi_1(x,\chi)$.
\begin{lemma}
Let $q\ge2$ be an integer, and let $x$, $\sigma$, and $T$ be real
numbers satisfying $\half\le\sigma<1$ and $1\le x\ll T\ll q$. If
$q\notin\Qs T$, then for all nonprincipal $\chi\in\Phi_*(q)$, we have
\begin{equation*}
\psip x\chi \ll x^{1+\sigma}\log q.
\end{equation*}
\label{psiplem}
\end{lemma}

\begin{proof}
We begin by writing
\begin{equation*}
\psi_1(x,\chi) = \frac{-1}{2\pi i} \int_{2-i\infty}^{2+i\infty} {L'\over
L}(s,\chi) {x^{s+1}ds\over s(s+1)}
\end{equation*}
and pulling the contour leftwards towards $\Re s=-\infty$ to see that
\begin{equation*}
\psip x\chi = -\sum_\rho {x^{\rho+1}\over\rho(\rho+1)} + O(x\log x),
\end{equation*}
where the sum runs over all nontrivial zeros $\rho=\beta+i\gamma$ of
$L(s,\chi)$ (see for instance \cite[Chapter 19]{Dav:MNT}). Because $q$
is not in $\Qs T$, every zero of $L(s,\chi)$ has either
$\beta\le\sigma$ or $\abs\gamma\ge T$, and thus we can write
\begin{equation*}
\psip x\chi \ll \sum_{\beta\le\sigma} {x^{1+\beta}\over\gamma^2} +
\sum_{\abs\gamma\ge T} {x^{1+\beta}\over\gamma^2} + x\log x.
\end{equation*}
However, the number of zeroes of $L(s,\chi)$ up to height $T$ is $\ll
T\log qT$, and so $\sum_{\abs\gamma\ge T} \gamma^{-2} \ll T^{-1}\log
qT$ by partial summation. Therefore
\begin{equation*}
\psip x\chi \ll x^{1+\sigma}\log q + x^2T^{-1}\log qT + x\log x.
\end{equation*}
Since $x\ll T\ll q$, the first term is dominant, and the lemma is
established.
\end{proof}

%% file: shifted.tex
Let $\A$ be a finite sequence, $\nu$ a map from $\A$ to the positive
integers, and $w$ a function from $\A$ to the nonnegative reals. Let
$\Upsilon$ be a squarefree integer, put
\begin{equation*}
S(\A,\Upsilon) = \sum \begin{Sb}a\in\A \\ (\nu(a),\Upsilon)=1\end{Sb}
w(a),
\end{equation*}
and, for all $d$ dividing $\Upsilon$, put
\begin{equation*}
A_d = \sum \begin{Sb}a\in\A \\ d\mid\nu(a)\end{Sb} w(a).
\end{equation*}
\begin{lemma}
Suppose that $X$ and $R$ are positive numbers and $f(d)$ a
multiplicative function such that for all $d$ dividing $\Upsilon$, we have
$f(d)\ge d$ and
\begin{equation}
\abs{A_d-\frac X{f(d)}} \le R.  \label{bddbyR}
\end{equation}
Then there exists an absolute positive constant $C_1$ such that
\begin{equation*}
S(\A,\Upsilon) \ge {C_1X\over\lg1\omega(\Upsilon)} \prod_{p\mid
\Upsilon} \big( 1-\frac1{f(p)} \big) + O\big( R(\omega(\Upsilon))^2 \big).
\end{equation*}
\label{shiftlem}
\end{lemma}

\begin{proof}
Let $p_j$ denote the $j$th prime, and put $z=p_{\omega(\Upsilon)}$ and
$P=\prod_{p\le z}p$. Also let $\{\lamm d\}$ be a sequence of real
numbers such that $\lamm 1\le1$ and, if we define
\begin{equation*}
\sigma_n = \sum_{d\mid n} \lamm d,
\end{equation*}
then $\sigma_n\le0$ for all integers $n\ge2$. We begin by citing the
lower bound
\begin{equation}
S(\A,\Upsilon) \ge X \prod_{p\mid \Upsilon} \big( 1-\frac1{f(p)} \big)
\sum_{d\mid P} {\sigma_d\over\prod_{p\mid d}(p-1)} - R \sum_{d\mid P}
\abs{\lamm d}.
\label{shiftbd}
\end{equation}
This is a special case of the shifted sieve of Iwaniec \cite[Lemma
1]{Iwa:OtPoJ}, where we have specified that $Q=\Upsilon$, $A=R$,
$B=1$, and $g(d)=d$ for all $d$ dividing $P$, and that the
correspondence $l$ sends the smallest prime factor of $\Upsilon$ to
$p_1$, the next smallest to $p_2$, and so on. We now take $\{\lamm
d\}$ to be Rosser's weights for the linear sieve, whose definition
depends on a positive parameter $y$ as follows. If $d$ is not
squarefree, define $\lamm d=0$. If $d=q_1\cdots q_r$ for primes
$q_1>\cdots>q_r$, define
\begin{equation*}
\lamm d = \begin{cases}
(-1)^r & \hbox{if $q_1\cdots q_{2l-1}q_{2l}^3<y$ for all $0\le
l\le r/2$,} \\
0 & \hbox{otherwise.}
\end{cases}
\end{equation*}
We will need the following facts about the sequence $\{\lamm d\}$
\cite[Lemma 2]{Iwa:OtPoJ}: if $4\le z^2\le y\le z^4$, then
\begin{equation*}
\sum_{d\mid P} \abs{\lamm d} \ll y(\log y)^{-2}
\end{equation*}
and
\begin{equation}
\sum_{d\mid P} {\sigma_d\over\prod_{p\mid d}(p-1)} = 2e^\gamma
\frac{\log(s-1)}s + O\big( {1\over\log y} \big),
\label{Oconst}
\end{equation}
where $s=(\log y)/(\log z)$. Applying this with $y=C_2z^2$ for $C_2$ a
positive constant gives us
\begin{equation}
2e^\gamma \frac{\log(s-1)}s + O\big( {1\over\log y} \big) =
{e^\gamma\log C_2\over\log z} \big( 1+O\big( {\log C_2\over\log z}
\big) \big) + O\big( {1\over\log z} \big) \ge {C_1\over\log z}
\label{withc1}
\end{equation}
for some positive constant $C_1$, if $C_2$ and $z$ are sufficiently
large. With these estimates, the lower bound (\ref{shiftbd}) becomes
\begin{equation*}
S(\A,\Upsilon) \ge {C_1X\over\log z} \prod_{p\mid \Upsilon} \big(
1-\frac1{f(p)} \big) + O\big( {RC_2z^2\over(\log z)^2} \big).
\end{equation*}
We note that $C_2$ is an absolute constant, since it depends only on
the $O$-constant in equation (\ref{Oconst}), and thus $C_1$ is
absolute as well, since it depends only on $C_2$ and the $O$-constants
in equation (\ref{withc1}). It remains only to note that $z\sim
\omega(\Upsilon)\lg1\omega(\Upsilon)$ to establish the lemma.
\end{proof}

%% file: applyshift.tex
We may now establish Proposition \ref{gqprop}. Let $q\ge2$ be an
integer and $x>1$ and $\half\le\sigma<1$ real numbers. We will apply
Lemma \ref{shiftlem} with $\A$ being the set of positive integers less
than $x$. Let $\Upsilon=s(\phi(q))$, let $\nu(n)$ be defined as in
Section \ref{charsec} before equation (\ref{manyBiggs2}), and let
$w(n)=\Lambda(n)(x-n)$. From the relation (\ref{manyBiggs2}), we see
that
\begin{equation*}
S(\A,\Upsilon) = \sum_{n<x} \gamma(n)\Lambda(n)(x-n)
\end{equation*}
counts only prime powers which are \lr{s}\mod q. Using the form
(\ref{charfnGd}) for $h_d$ and the definition of the $\psip x\chi$, we
also have
\begin{equation}
\begin{split}
A_d &= \sum \begin{Sb}n<x \\ d\mid \nu(n)\end{Sb} w(n) = \sum_{n<x}
h_d(n)w(n) \\
&= \frac1{\abs{\Phi_d}} \sum_{\chi\in\Phi_d} \sum_{n<x} \chi(n)w(n) \\
&= \frac1{d^{m(d)}} \psip x{\chi_0} + \frac1{\abs{\Phi_d}}
\sum \begin{Sb}\chi\in\Phi_d \\ \chi\ne\chi_0\end{Sb} \psip x\chi.
\end{split}
\label{apply1}
\end{equation}
If we write $\psi_1(x) = \sum_{n<x} \Lambda(n)(x-n)$, then
\begin{equation*}
\psi_1(x) - \psip x{\chi_0} = \sum \begin{Sb}n<x \\
(n,q)>1\end{Sb} \Lambda(n)(x-n) \ll x \sum_{p\mid q} \sum
\begin{Sb}r\ge1 \\ p^r<x\end{Sb} \log p \ll (x\log x)\log q,
\end{equation*}
since $\omega(q)\ll\log q$. Moreover, if we assume that $q\notin\Qs
x$, then we may apply Lemma \ref{psiplem} (with $T=x$) to bound the
terms in the last sum of equation (\ref{apply1}); we obtain
\begin{equation*}
A_d = \frac1{d^{m(d)}} \psi_1(x) + O\big( x^{1+\sigma}\log q \big).
\end{equation*}
Thus if we take $X=\psi_1(x)$ and $f(d)=d^{m(d)}$ for all $d$ dividing
$s(\phi(q))$, we see that we can take $R\ll x^{1+\sigma}\log
q$. Applying Lemma \ref{shiftlem}, we see that
\begin{equation*}
\begin{split}
S(\A,\Upsilon) &\ge {C_1\psi_1(x)\over\lg1\w}c_0 + O\big(
(x^{1+\sigma}\log q) \w^2 \big) \\
&= {C_1\psi_1(x)\over\lg1\w}c_0 \big( 1+O\big( x^{-1+\sigma} \big( \w^2\lg1\w
\big)c_0^{-1}\log q \big) \big) \\
&= {C_1\psi_1(x)\over\lg1\w}c_0 \big( 1+O\big(\big( x^{-1}f(q,\sigma)
\big)^{1-\sigma}\big) \big),
\end{split}
\end{equation*}
since the bound $\psi_1(x)\gg x^2$ follows from Chebyshev's bound for
$\psi(x)$. Assuming that $x$ exceeds a sufficiently large (in terms of
$\sigma$) multiple of $f(q,\sigma)$, we obtain a positive lower bound
for $S(\A,\Upsilon)$. Therefore, there exists a prime power
$p^r\ll_\sigma f(q,\sigma)$ which is a \lr\mod q. But if $p^r$ is a
\lr, we must have $(r,\phi(q))=1$, in which case $p$ itself is also a
\lr{} which is $\ll_\sigma f(q,\sigma)$. This establishes the
proposition.

%% file: Zhang.tex
To establish Lemma \ref{aaplem}, we introduce the notation $\Qps T$ to
denote the subset of $\Qs T$ consisting of the odd prime powers, and
we recall that $Q'(Y;\sigma,T)$ denotes the number of elements of
$\Qps T$ not exceeding $Y$. Given an odd prime power $p^r$, every
character in $\Phi_*(p^r)$ is induced by a character\mod{p^2}
\cite[Lemma 6]{Mar:UBftLA-PPR}. The proof of this fact is similar to
the proof that any primitive root\mod{p^2} is also a primitive
root\mod{p^r} for every odd prime $p$ and integer $r\ge3$.

Consequently, for every prime power $p^r\in\Qps T$, there is a
character $\chi$ which is primitive to one of the moduli $p$ or $p^2$
such that $L(s,\chi)$ has a zero $\beta+i\gamma$ with $\beta>\sigma$
and $\abs\gamma<T$. On the other hand, every such character will
account for $\ll\log Y$ prime powers in $\Qps T$ which do not exceed
$Y$, and so
\begin{equation}
Q'(Y;\sigma,T) \ll (\log Y) \sum_{q<Y} \starsum_{\chi\mod q}
N(\sigma,T,\chi),  \label{logY}
\end{equation}
where $N(\sigma,T,\chi)$ denotes the number of zeros $\beta+i\gamma$
of $L(s,\chi)$ satisfying $\beta>\sigma$ and $\abs{\gamma}<T$, and
$\starsum$ denotes a summation over primitive characters only. Zhang
\cite{Zha:OtDoZotDL-f} has established the following zero-density
estimate for Dirichlet $L$-functions: for any real numbers $Y$,
$\delta>0$ and $\frac{17}{22}\le\sigma\le1$, we have
\begin{equation}
\sum_{q<Y} \starsum_{\chi\mod q} N(\sigma,T,\chi) \ll_\delta
(Y^2T)^{6(1-\sigma)/(5\sigma-1)+\delta}.
\label{Zhangthm}
\end{equation}
We apply this estimate with $T=Y^{\theta}$ and $\sigma=1-B^{-1}$,
where $B$ is as in Theorem \ref{aapthm}. Together with the bound
(\ref{logY}), this gives us $Q'(Y;\sigma,T) \ll_{\ep,\eta} Y^{\ep}$,
as long as $\delta=\delta(\ep,\eta)$ and $\theta=\theta(\ep,\eta)$ are
small enough with respect to $\ep$ and $\eta$. This establishes Lemma
\ref{aaplem}.

%% file: almostall.tex
Unfortunately, a given character can in general induce characters in
$\Phi_*(q)$ for many more moduli $q$ if we do not restrict to prime
powers, and so we must work harder to establish Lemma \ref{aalem}.
Given positive integers $m$ and $n$ such that $m$ divides $n$, we say
that $n$ is an {\it admissible multiple} of $m$ if there exists a
character in $\Phi_*(n)$ which is induced by a primitive character\mod
m.

\begin{lemma}
Let $q\ge2$ be an integer, and set $t=\omega(q)$. Let
$p_1$,\ldots\!\!, $p_t$ be the primes dividing $q$ and
$r_1$,\ldots\!\!, $r_t$ positive integers. Then for every admissible
multiple $nq$ of $q$, either:
\begin{equation*}
\begin{split}
(i)& \hbox{ $p_i^{r_i}$ divides $n$ for some $1\le i\le t$; or} \\
(ii)& \hbox{ $n$ is not divisible by any prime which is congruent to
$1\mod{\phi^2(q)p_1^{r_1}\cdots p_t^{r_t}}$}.
\end{split}
\end{equation*}
\label{craftylem}
\end{lemma}

\def\chim#1#2{{\chi_{#1}^{(#2)}}}
\begin{proof}
We use parenthetical superscripts to indicate explicitly the modulus of
a character, so that $\chim{}q$ denotes a character\mod q, for
example. To establish the lemma, it suffices to show that if (i) and
(ii) both fail, then any character $\chim{}q$ which induces an element
$\chim1{nq}$ of $\Phi_*(nq)$ is in fact principal (hence imprimitive),
contradicting the assumption that $nq$ is an admissible multiple of
$q$.

Assume the negations of (i) and (ii). Write $nq=n'q'$, where $q'$ is
the largest divisor of $nq$ with $s(q')=s(q)$, so that $q$ divides
$q'$ and $(n',q')=1$. Then any character\mod{nq} is the product of a
character\mod{n'} and a character\mod{q'}. Since
$\chim1{nq}\in\Phi_*(nq)$, we may write
\begin{equation*}
\chim1{nq}=\big( \chim2{n'}\chim3{q'} \big)^{E(nq)/s(E(nq))}
\end{equation*}
for some characters $\chim2{n'}$ and $\chim3{q'}$. Since $p_i^{r_i}$
does not divide $n$ for any $1\le i\le t$, we see from the definition
of $q'$ that $\phi(q')$ divides $\phi(q)p_1^{r_1-1}\cdots
p_t^{r_t-1}$. On the other hand, $n$ is divisible by a prime which is
congruent to 1\mod{\phi^2(q)p_1^{r_1}\cdots p_t^{r_t}}, and so
$\phi^2(q)p_1^{r_1}\cdots p_t^{r_t}$ must divide $E(nq)$. These
observations together imply that $\phi(q')$ divides $E(nq)/s(E(nq))$,
and thus
\begin{equation*}
\big( \chim2{n'}\chim3{q'} \big)^{E(nq)/s(E(nq))} = \big( \chim2{n'}
\big)^{E(nq)/s(E(nq))}\chim0{q'},
\end{equation*}
where $\chim0{q'}$ is the principal character\mod{q'}. We see that the
character $\chim1{nq}$ induced by $\chim{}q$ is also induced by a
character\mod{n'}. But since $(q,n')=1$, it must be the case that
$\chim{}q$ is principal. This establishes the lemma.
\end{proof}

Let $A(x;q)$ be the number of admissible multiples of $q$ not
exceeding $x$.
\begin{lemma}
Let $\delta>0$ be a real number and $x$, $y=y(x)$, and $z=z(x)$ real
parameters satisfying $x$, $y$, $z>1$ and
\begin{equation}
z^3y^{\log z} \ll (\log x)^{1-\delta}.  \label{uniform}
\end{equation}
Then for all integers $q$ with $2\le q\le z$, we have
\begin{equation}
A(xq;q) \ll_\delta \frac{x\log z}y + \frac x{\exp( \lg2 x/z^3y^{\log
z} )}.  \label{yzbound}
\end{equation}
\label{liladmiss}
\end{lemma}

\begin{proof}
Set $t=\omega(q)$, and choose integers $r_i$ such that
\begin{equation}
p_i^{r_i-1} \le y \le p_i^{r_i} \quad (1\le i\le t). \label{ridef}
\end{equation}
By applying Lemma \ref{craftylem}, we see that the number of
admissible multiples $nq$ of $q$ with $n<x$ is bounded by
\begin{equation}
\sum_{i=1}^t \frac x{p_i^{r_i}} + \#\{ n<x : p\mid n\implies
p\not\equiv1\mod{\phi^2(q)p_1^{r_1}\cdots p_t^{r_t}} \}.
\label{twosums}
\end{equation}
In the first term, we use the estimate $t\le\log z$ for $z$
sufficiently large, and the choice (\ref{ridef}) of the $r_i$, to see
that
\begin{equation}
\sum_{i=1}^t \frac x{p_i^{r_i}} \le \frac{x\log z}y. \label{firstsum}
\end{equation}
We treat the second term using a simple upper bound sieve. Notice that
by the choice (\ref{ridef}) of the $r_i$, we have
\begin{equation}
\phi^2(q)p_1^{r_1}\cdots p_t^{r_t} \le q^2 \bigg( \prod_{i=1}^t yp_i \bigg)
\le q^2(y^tz) \le z^3y^{\log z}.
\label{modsmall1}
\end{equation}
The prime number theorem for arithmetic progressions states that
given $\delta>0$, we have
\begin{equation*}
\psi(x;d,1) = \frac x{\phi(d)} + O_\delta\big( x\exp( -C_3(\log x)^{1/2}) \big)
\end{equation*}
for some positive constant $C_3$, uniformly for all $d\ll(\log
x)^{1-\delta}$ \cite[equations (10)--(11) of Section 20]{Dav:MNT}. By
partial summation, this implies that
\begin{equation}
\sum \begin{Sb}p<x \\ p\equiv1\mod d\end{Sb} p^{-1} = {\lg2
x\over\phi(d)} + O_\delta(1),
\label{sumpinv}
\end{equation}
again uniformly for $d$ in the above range, which includes
$d=\phi^2(q)p_1^{r_1}\cdots p_t^{r_t}$ due to equation
(\ref{modsmall1}) and the restriction (\ref{uniform}). The formula
(\ref{sumpinv}) allows us to apply an upper bound sieve from
Halberstam--Richert \cite[Corollary 2.3.1]{HalRic:SM} to deduce that
\begin{equation*}
\#\{ n<x : p\mid n\implies p\not\equiv1\mod{\phi^2(q)p_1^{r_1}\cdots
p_t^{r_t}} \} \ll_\delta x(\log x)^{-1/\phi(\phi^2(q)p_1^{r_1}\cdots
p_t^{r_t})}.
\end{equation*}
We rewrite this using the bound (\ref{modsmall1}) as
\begin{equation*}
\#\{ n<x : p\mid n\implies p\not\equiv1\mod{\phi^2(q)p_1^{r_1}\cdots
p_t^{r_t}} \} \ll_\delta \frac x{\exp( \lg2 x/z^3y^{\log z})}.
\end{equation*}
Using this bound together with the bound (\ref{firstsum}) in equation
(\ref{twosums}) establishes the lemma.
\end{proof}

\def\back{\hskip-.5cm}
Define $\R \sigma T$ to be the set of integers $q\ge3$ such that, for
some primitive character $\chi\mod q$, the corresponding $L$-function
$L(s,\chi)$ has a zero $\beta+i\gamma$ with $\beta>\sigma$ and
$\abs\gamma<T$.
\begin{lemma}
For all real $x>1$, we have
\begin{equation}
\sum \begin{Sb}q<x \\ q\in\R{17/22}{x^{1/20}}\end{Sb}\back 1 \ll
x^{.997} \quad\hbox{and}\quad \back\!\sum \begin{Sb}x<q \\
q\in\R{17/22}{x^{1/20}}\end{Sb}\back \!q^{-1} \ll x^{-.003}.
\label{RsTbd}
\end{equation}
\label{Zhanglem}
\end{lemma}

\begin{proof}
The right-hand side of the zero-density estimate (\ref{Zhangthm}) is
certainly an upper bound for the first sum in (\ref{RsTbd}) as
well. Taking $Y=x$, $T=x^{1/20}$, and $\theta=\frac1{100}$ in
(\ref{Zhangthm}), we see that
\begin{equation*}
\sum \begin{Sb}q<x \\ q\in\R{17/22}{x^{1/20}}\end{Sb}\back 1
\ll x^{41861/42000},
\end{equation*}
and $41861/42000<.997$. This establishes the first bound in
(\ref{RsTbd}), and the second bound follows directly by partial
summation.
\end{proof}

We are now ready to prove Lemma \ref{aalem}. We note that every
element of $\Qs T$ is an admissible multiple of some element of
$\R\sigma T$. Therefore,
\begin{equation}
Q(Y;\sigma,T) \le \sum \begin{Sb}q<Y \\ q\in\R\sigma T\end{Sb} A(Y;q).
\label{usingadmiss}
\end{equation}
For $q\le\lg3Y$, we bound $A(Y;q)$ by applying Lemma \ref{liladmiss}
with $z=\lg3Y$ and $y=(\lg2Y)^{1/(2\log z)}$, which satisfy the
condition (\ref{uniform}) with any $\delta<1$. Of the two terms
in equation (\ref{yzbound}), the first term is dominant, giving
\begin{equation*}
A(Y;q) \le A(Yq;q) \ll \frac{Y\lg4 Y}{\exp( \lg3Y/2\lg4Y )}.
\end{equation*}
For the remaining values of $q$, we have the trivial bound $A(Y;q)\le
Y/q$. Therefore equation (\ref{usingadmiss}) becomes
\begin{equation*}
Q(Y;\sigma,T) \ll \sum_{q<\lg3 Y} \frac{Y\lg4 Y}{\exp(
\lg3Y/2\lg4Y )} + \sum \begin{Sb}\lg3 Y\le q<Y \\ q\in\R\sigma
T\end{Sb} \frac Yq.
\end{equation*}
Upon choosing $\sigma=\frac{17}{22}$ and $T=Y^{1/20}$, we apply Lemma
\ref{Zhanglem} to the second sum to obtain
\begin{equation*}
Q\big( Y; \hbox{$17\over22$} ,Y^{1/20} \big) \ll \frac{Y\lg3Y\lg4Y}
{\exp( \lg3Y/2\lg4Y )} + \frac Y{(\lg3Y)^{.003}} = o(Y),
\end{equation*}
which establishes the lemma.